\documentclass[12pt,leqno]{article}
\usepackage{amsmath, amssymb}

\newcommand{\sect}[1]{\setcounter{equation}{0}\section{#1}}

\def\epsilon{\varepsilon}

\begin{document}

\LARGE \noindent 
{\bf Cartan-Thullen theorem \\ 
for a $\mathbb C^n$-holomorphic function
\\ and a related problem}

\large

\vspace*{0.4em}

\hfill Hiroki Yagisita (Kyoto Sangyo University)

\vspace*{1.2em}

\normalsize

Abstract:

Cartan-Thullen theorem is a basic one in the theory 
of analytic functions of several complex variables. 
It states that for any open set $U$ of ${\mathbb C}^k$, 
the following conditions are equivalent:
(a) $U$ is a domain of existence,
(b) $U$ is a domain of holomorphy
and (c) $U$ is holomorphically convex.
On the other hand, when $f \, (\, =(f_1,f_2,\cdots,f_n)\, )$ 
is a $\mathbb C^n$-valued function 
on an open set $U$ 
of $\mathbb C^{k_1}\times\mathbb C^{k_2}\times\cdots\times\mathbb C^{k_n}$, 
$f$ is said to be $\mathbb C^n$-analytic, 
if $f$ is complex analytic and for any $i$ and $j$, 
$i\not=j$ implies $\frac{\partial f_i}{\partial z_j}=0$, 
where $(z_1,z_2,\cdots,z_n)
\in \mathbb C^{k_1}\times\mathbb C^{k_2}\times\cdots\times\mathbb C^{k_n}$ holds. 
We note that  a $\mathbb C^n$-analytic mapping 
and a $\mathbb C^n$-analytic manifold can also be easily defined. 

In this paper, we show an analogue of Cartan-Thullen theorem
for $\mathbb C^n$-analytic functions.  
For $n=1$, it gives Cartan-Thullen theorem itself. 
Our proof is almost the same as Cartan-Thullen theorem. 
Thus, our generalization seems to be natural. 
On the other hand, our result is partial, 
because we do not answer the following question. 
That is, does a connected open $\mathbb C^n$-holomorphically convex set $U$ 
exist such that $U$ is not the direct product 
of any holomorphically convex sets 
$U_1, U_2, \cdots, U_{n-1}$ and $U_n$ ? 
As a corollary of our generalization, we give the following partial result. 
If $U$ is convex, then $U$ is the direct product 
of some holomorphically convex sets. 

Also, $f$ is said to be $\mathbb C^n$-triangular, 
if $f$ is complex analytic and for any $i$ and $j$, 
$i<j$ implies $\frac{\partial f_i}{\partial z_j}=0$. 
Kasuya suggested that a $\mathbb C^n$-analytic manifold 
and a $\mathbb C^n$-triangular manifold might, for example, 
be related to a holomorphic web and a holomorphic foliation.  

\vspace*{0.8em} 
\noindent 
Keywords: \ 

\noindent 
Stein space, pseudoconvex manifold, holomorphic foliation, holomorphic web. 

\vfill 
\noindent 

\vspace*{0.4em} 

\noindent 

\newpage

\sect{Introduction}
\ \ \ \ \ 
First, we generalize the notion of a holomorphic function. 

\noindent 
{\bf Definition 1} (Structure sheaf) : 

Let $k_1, k_2, \cdots, k_{n-1}, k_n, l_1, l_2, \cdots, l_{n-1}$ 
and $l_n$ be natural numbers. 
Let $U$ be an open set of 
$\mathbb C^{l_1}\times\mathbb C^{l_2}\times\cdots\times\mathbb C^{l_n}$. 
Let $f \, (\, =(f_1,f_2,\cdots,f_n)\, )$ be a map from $U$ 
to $\mathbb C^{k_1}\times\mathbb C^{k_2}\times\cdots\times\mathbb C^{k_n}$. 
Then, $f$ is said to be $\mathbb C^n$-holomorphic ($\mathbb C^n$-analytic), 
if $f$ is holomorphic and for any $a\in U$ and any $i, j \in \{1,2,\cdots,n\}$, 
$i\not=j$ implies $\frac{\partial f_i}{\partial z_j}(a)=0$, 
where $(z_1,z_2,\cdots,z_n)
\in \mathbb C^{l_1}\times\mathbb C^{l_2}\times\cdots\times\mathbb C^{l_n}$ holds. 

Let $O_{l_1,l_2,\cdots,l_n}(U)$ denote the set 
of all $\mathbb C^n$-valued $\mathbb C^n$-holomorphic functions on $U$.  
Then, $\{O_{l_1,l_2,\cdots,l_n}(U)\}_U$ is called the sheaf of germs 
of $\mathbb C^n$-holomorphic functions. 

\noindent 
{\bf Example 2} : \ 

(1) \ Let $$\pi_j(U):=\{z_j \in \mathbb C^{l_j}
|\exists z_1,z_2,\cdots,z_{j-1}, z_{j+1},z_{j+2},\cdots,z_n: 
(z_1,z_2,\cdots,z_n)\in U\}.$$
Let $f_j$ be a holomorphic function on $\pi_j(U)$. 
Then, $(f_1,f_2,\cdots,f_n)$ 
is a $\mathbb C^n$-holomorphic function on $U$. 

(2) \ Let $\varepsilon$ be a small positive number. 
Let $$U \, := \, \cup_{\theta \in \mathbb R} \, 
(\, \{ \, z_1\in \mathbb C \, | \, |z_1-e^{\sqrt{-1}\theta}|<\varepsilon\, \}
\, \times\, \{\, z_2\in \mathbb C\, |\, |z_2-\theta|<\varepsilon\, \}\, ).$$
Then, $(\log z_1,0)$ is a $\mathbb C^2$-holomorphic function on $U$. 
However, $\log z_1$ is a multivalued function on $\pi_1(U)$.

\noindent 
{\bf Remark 3} : \ 

(1) \ The composition of $\mathbb C^n$-holomorphic mappings 
is $\mathbb C^n$-holomorphic. So, a $\mathbb C^n$-analytic manifold 
can be easily defined with its structure sheaf.

(2) \  For $n=1$, $\{O_{l}(U)\}_U$ is the sheaf of germs 
of holomorphic functions. 

(3) \ 
$(f_1,f_2,\cdots,f_n)$ is $\mathbb C^n$-holomorphic, 
if and only if $(f_1,0,0,\cdots,0,0,0)$, 
$(0,f_2,0,\cdots,0,0,0)$, $\cdots$, 
$(0,0,0,\cdots,0,f_{n-1},0)$ 
and $(0,0,0,\cdots,0,0,f_n)$ 
are $\mathbb C^n$-holomorphic. Also, 
$(f_1,f_2,\cdots,f_n), (g_1,g_2,\cdots,g_n) \in O_{l_1,l_2,\cdots,l_n}(U)$ 
implies 
$(f_1g_1,f_2g_2,\cdots,f_ng_n) \in O_{l_1,l_2,\cdots,l_n}(U)$. 
Further, if a sequence $\{f_m\}_{m=1}^\infty$ in $O_{l_1,l_2,\cdots,l_n}(U)$ 
uniformly convergences to $g\in (O_{l_1+l_2+\cdots+l_n}(U))^n$ on compact sets, 
then $g\in O_{l_1,l_2,\cdots,l_n}(U)$ holds. 
So, $O_{l_1,l_2,\cdots,l_n}(U)$ is the direct product 
of closed $\mathbb C$-subalgebras of the usual one $O_{l_1+l_2+\cdots+l_n}(U)$.

(4) \ 
When $A$ is a commutative Banach algebra, Lorch ([6]) 
gave a definition that an $A$-valued function on an open set of $A$
is $A$-holomorphic. With the norm $\max_{j=1,2,\cdots,n}|z_j|$, 
$\mathbb C^n$ is a locally compact one.  
We did a little study on $A$-analytic manifolds ([15, 16]). 
\hfill ---

\newpage 

Since the structure sheaf $\{O_{l_1,l_2,\cdots,l_n}(U)\}_U$ 
was defined, we define $\mathbb C^n$-existence, 
$\mathbb C^n$-holomorphy
and $\mathbb C^n$-holomorphic convexity. 
Just in case, we state uniqueness theorem. 

\noindent 
{\bf Proposition 4} : 

Let $U$ be a connected open set of 
$\mathbb C^{l_1}\times\mathbb C^{l_2}\times\cdots\times\mathbb C^{l_n}$. 
Let $f, g \in O_{l_1,l_2,\cdots,l_n}(U)$. 
Let $a\in U$. If for any multi-index $\alpha$, 
$\frac{\partial^{|\alpha|} f}{{\partial z}^\alpha}(a)
=\frac{\partial^{|\alpha|} g}{{\partial z}^\alpha}(a)$ holds, 
then $f=g$ holds. 

\noindent 
{\sf Proof} : \ It is an easy corollary of the usual uniqueness theorem. 
\hfill 
$\blacksquare$

\noindent 
{\bf Definition 5} (Existence, Holomorphy) : 

Let $U$ be an open set of 
$\mathbb C^{l_1}\times\mathbb C^{l_2}\times\cdots\times\mathbb C^{l_n}$. 

(1) \ $U$ is said to be a domain of $\mathbb C^n$-existence, 
if the following holds. 
There exists $f\in O_{l_1,l_2,\cdots,l_n}(U)$ 
such that for any open sets $V$ and $W$ of
$\mathbb C^{l_1}\times\mathbb C^{l_2}\times\cdots\times\mathbb C^{l_n}$, 
if $V$ is connected and $\emptyset\not=V\setminus U$ and 
$\emptyset\not=W\subset U\cap V$ hold, 
then for any $g\in  O_{l_1,l_2,\cdots,l_n}(V)$, 
$f_{\upharpoonright W}\not=g_{\upharpoonright W}$ holds.  

(2) \ $U$ is said to be a domain of $\mathbb C^n$-holomorphy, 
if the following holds. For any open sets $V$ and $W$ of
$\mathbb C^{l_1}\times\mathbb C^{l_2}\times\cdots\times\mathbb C^{l_n}$, 
if $V$ is connected and $\emptyset\not=V\setminus U$ and 
$\emptyset\not=W\subset U\cap V$ hold, 
then there exists $f\in O_{l_1,l_2,\cdots,l_n}(U)$ 
such that for any $g\in  O_{l_1,l_2,\cdots,l_n}(V)$, 
$f_{\upharpoonright W}\not=g_{\upharpoonright W}$ holds.  

\noindent 
{\bf Lemma 6} : 

$\mathbb C^n$-existence implies $\mathbb C^n$-holomorphy. 

\noindent 
{\sf Proof} : \ It is obvious. 
\hfill 
$\blacksquare$

\noindent 
{\bf Definition 7} (Holomorphic convexity) : 

Let $|\{w_k\}_{k=1}^m|$ denote $\max_{k=1,2,\cdots,m}|w_k|$ 
for $w_1,w_2,\cdots,w_m \in \mathbb C$. 
Let $U$ be an open set of 
$\mathbb C^{l_1}\times\mathbb C^{l_2}\times\cdots\times\mathbb C^{l_n}$.

(1) \ Let $K$ be a compact subset of $U$. 
Let 
$${\widehat K}_{l_1,l_2,\cdots,l_n}^U \, :=\, \{ \, z\in U \, 
| \, \forall \, f \, \in \, O_{l_1,l_2,\cdots,l_n}(U): \, 
|f(z)|\leq \sup_{w\in K}|f(w)| \, \}.$$
Then, ${\widehat K}_{l_1,l_2,\cdots,l_n}^U$ 
is called the $\mathbb C^n$-holomorphically convex hull of $K$. 

(2) \ $U$ is said to be $\mathbb C^n$-holomorphically convex, 
if for any compact subset $K$ of $U$, 
${\widehat K}_{l_1,l_2,\cdots,l_n}^U$ is compact. 
\hfill ---

\vspace*{0.8em} 

The following is the main result. 
We note that for $n=1$, it is Cartan-Thullen theorem ([1]) itself. 

\noindent 
{\bf Theorem 8} : 

Let $U$ be an open set of 
$\mathbb C^{l_1}\times\mathbb C^{l_2}\times\cdots\times\mathbb C^{l_n}$. 
Then, the following conditions are equivalent: 
(a) $U$ is a domain of $\mathbb C^n$-existence,
(b) $U$ is a domain of $\mathbb C^n$-holomorphy
and (c) $U$ is $\mathbb C^n$-holomorphically convex.

\noindent 
{\bf Remark 9} : 

Let $U_j \, (\not=\emptyset)$ 
be a connected open set of $\mathbb C^{l_j}$ $(j=1,2,\cdots,n)$. 
Let $U:=U_1\times U_2\times \cdots \times U_n$. 

(1) \ Let $K_j$ be a compact subset of $U_j$ $(j=1,2,\cdots,n)$. 
Then, $$\widehat{(K_1\times K_2 \times \cdots \times K_n)}_{l_1,l_2,\cdots,l_n}^U
\, = \, \widehat{K_1}_{l_1}^{U_1}\times
\widehat{K_2}_{l_2}^{U_2}\times \cdots 
\times\widehat{K_n}_{l_n}^{U_n}$$ 
holds. 

(2) \ $U$ is $\mathbb C^n$-holomorphically convex, 
if and only if $U_1, U_2, \cdots, U_{n-1}$ and $U_n$ 
are holomorphically convex. 

\noindent 
{\sf Proof} : \ (1) \ $U_1\times U_2\times \cdots \times U_{j-1}
\times U_{j+1}\times U_{j+2}\times\cdots\times U_n$ 
is connected. Hence, if $(f_1,f_2,\cdots,f_n) \in O_{l_1,l_2,\cdots,l_n}(U)$ 
holds, then for any $a_j\in U_j$, the function 
$(z_1,z_2,\cdots,z_{j-1},z_{j+1},z_{j+2},\cdots,z_n) 
\, \mapsto \, f_j(z_1,z_2,\cdots,z_{j-1},a_j,z_{j+1},z_{j+2},\cdots,z_n)$
is constant. So, $O_{l_1,l_2,\cdots,l_n}(U) \, = \, 
O_{l_1}(U_1)\times O_{l_2}(U_2)\times\cdots\times O_{l_n}(U_n)$ holds. 
For any $(z_1,z_2,\cdots,z_n)\in U$, 
$$\forall \, f\in O_{l_1,l_2,\cdots,l_n}(U):\, |f(z_1,z_2,\cdots,z_n)| \, \leq \, 
\sup_{w\in K_1\times K_2 \times \cdots \times K_n}|f(w)|$$
$$\Longleftrightarrow$$ 
$$\forall \, (f_1,f_2,\cdots,f_n)\in 
O_{l_1}(U_1)\times O_{l_2}(U_2)\times\cdots\times O_{l_n}(U_n)$$
$$:\, \max_{i=1,2,\cdots,n} |f_i(z_i)| \, \leq \, 
\max_{i=1,2,\cdots,n} (\sup_{w_i\in K_i}|f_i(w_i)|)$$
$$\Longleftrightarrow$$ 
$$\forall \, i\in\{1,2,\cdots,n\}, \, \forall \, f_i\in O_{l_i}(U_i)
:\, |f_i(z_i)| \, \leq \, \sup_{w_i\in K_i}|f_i(w_i)|$$
holds. 

\ (2) \  Suppose that $U$ is $\mathbb C^n$-holomorphically convex. 
We show that $U_j$ is holomorphically convex. 
Let $K_j$ be a compact subset of $U_j$. 
There exists $(a_1,a_2,\cdots,a_n)\in U$. 
From (1), 
$$\widehat{(\{a_1\}\times \{a_2\} \times \cdots \times \{a_{j-1}\} \times K_j \times 
\{a_{j+1}\} \times \{a_{j+2}\} \times \cdots \times \{a_n\})}_{l_1,l_2,\cdots,l_n}^U$$ 
$$= \, \widehat{\{a_1\}}_{l_1}^{U_1}\times
\widehat{\{a_2\}}_{l_2}^{U_2}\times \cdots \times 
\widehat{\{a_{j-1}\}}_{l_{j-1}}^{U_{j-1}}\times 
\widehat{K_j}_{l_j}^{U_j}\times 
\widehat{\{a_{j+1}\}}_{l_{j+1}}^{U_{j+1}}\times
\widehat{\{a_{j+2}\}}_{l_{j+2}}^{U_{j+2}}\times \cdots \times 
\widehat{\{a_n\}}_{l_n}^{U_n}$$ 
holds. Hence, 
$$\pi_j(\widehat{(\{a_1\}\times \{a_2\} \times \cdots \times \{a_{j-1}\} \times K_j \times 
\{a_{j+1}\} \times \{a_{j+2}\} \times \cdots \times \{a_n\})}_{l_1,l_2,\cdots,l_n}^U) 
\, = \, \widehat{K_j}_{l_j}^{U_j}$$ 
holds. Because $U$ is $\mathbb C^n$-holomorphically convex, 
$\widehat{K_j}_{l_j}^{U_j}$ is compact. $U_j$ is holomorphically convex. 

Suppose that $U_1, U_2, \cdots, U_{n-1}$ and $U_n$ 
are holomorphically convex. We show that $U$ 
is $\mathbb C^n$-holomorphically convex. 
Let $K$ be a compact subset of $U$. Then, 
there exists $\{K_j\}_{j=1}^n$ such that $K_j$ 
is a compact subset of $U_j$ 
and $K\subset K_1\times K_2 \times \cdots \times K_n$ holds. 
So, from (1), 
$$\widehat{K}_{l_1,l_2,\cdots,l_n}^U 
\, \subset \, \widehat{K_1}_{l_1}^{U_1}\times
\widehat{K_2}_{l_2}^{U_2}\times \cdots 
\times\widehat{K_n}_{l_n}^{U_n}
\, ( \, \subset \, U \, )$$ 
holds. Because $U_1, U_2, \cdots, U_{n-1}$ and $U_n$ 
are holomorphically convex, $\widehat{K}_{l_1,l_2,\cdots,l_n}^U$ 
is compact. $U$ is $\mathbb C^n$-holomorphically convex. 
\hfill 
$\blacksquare$

Our generalization is considered natural. 
On the other hand, our result is partial, 
because we do not answer the following question. 

\noindent 
{\bf Question} : \ 

Does a connected $\mathbb C^n$-holomorphically convex open set 
(or, manifold) $U$ 
exist such that $U$ is not the direct product 
of any holomorphically convex ones 
$U_1, U_2, \cdots, U_{n-1}$ and $U_n$ ? 
\hfill ---

Now, we can give the following partial result. 

\noindent 
{\bf Corollary 10} : 

Let $U$ be a convex open set 
of $\mathbb C^{l_1}\times\mathbb C^{l_2}\times\cdots\times\mathbb C^{l_n}$. 

(1) \ Let $f\in O_{l_1,l_2,\cdots,l_n}(U)$. Then, 
there exists $g\in O_{l_1,l_2,\cdots,l_n}(\pi_1(U)\times \pi_2(U)\times \cdots \times \pi_n(U))$ 
such that $f=g_{\upharpoonright U}$ holds. 

(2) \  Suppose that $U$ is $\mathbb C^n$-holomorphically convex.  
Then, $U\, =\, \pi_1(U)\times \pi_2(U)\times \cdots \times \pi_n(U)$ holds. 
 
\noindent 
{\sf Proof} : \ (1) \ Let $f=(f_1,f_2,\cdots,f_n)$.  
For any $a_j\in \pi_j(U)$, $U\cap\pi_j^{-1}(\{a_j\})$
is convex, so, it is connected and the function 
$$(z_1,z_2,\cdots,z_{j-1},z_{j+1},z_{j+2},\cdots,z_n) 
\, \in \, U\cap\pi_j^{-1}(\{a_j\})$$
$$\mapsto \, \, \, f_j(z_1,z_2,\cdots,z_{j-1},a_j,z_{j+1},z_{j+2},\cdots,z_n) 
\, \in \, \mathbb C$$
is constant.  From this, it follows. 

(2) \ From Theorem 8, $U$ is a domain of $\mathbb C^n$-existence. 
Hence, from (1), it follows. 
\hfill 
$\blacksquare$

\vspace*{0.8em} 

\noindent 
{\bf Comment} : \ 

A map $f$ is said to be $\mathbb C^n$-triangular, 
if $f$ is holomorphic and for any $i$ and $j$, 
$i<j$ implies $\frac{\partial f_i}{\partial z_j}=0$. 
Kasuya suggested that a $\mathbb C^n$-analytic manifold 
and a $\mathbb C^n$-triangular manifold might, for example, 
be related to a holomorphic web and a holomorphic foliation.  
\hfill ---

\newpage 

\sect{Proof of main result} 
\ \ \ \ \ 
The proof of Theorem 8 is almost the same as Cartan-Thullen theorem. 
Perhaps, it seems to be also proved as a consequence of some general theory. 
However, for the sake of confirmation, we describe it. 
That is, we choose a proof that works in our case. 
In fact, it is extremely easy as we see below. 
When a reader believes that some proof which he knows works,  
he should skip the following proof. 

\vspace*{0.8em}

\noindent 
{\bf Lemma 11} : 

Let $K$ be a compact subset of $U$. 
Then, ${\widehat K}_{l_1,l_2,\cdots,l_n}^U$ is bounded. 

\noindent 
{\sf Proof} : \ Let $1\leq k \leq l_j$. Then, 
$(0,0,\cdots,0,z_{j,k},0,0,\cdots,0)\in O_{l_1,l_2,\cdots,l_n}(U)$ 
holds. Here, $z_j=(z_{j,1},z_{j,2},\cdots,z_{j,l_j})$ holds. 
Hence, $z\in {\widehat K}_{l_1,l_2,\cdots,l_n}^U$ 
implies $|z_{j,k}| \, \leq \, \sup_{w\in K}|w_{j,k}| \, (\, < \, +\infty\, )$. 
\hfill 
$\blacksquare$
 
\noindent 
{\bf Lemma 12} : 

Let $K$ be a compact subset of $U$. 
Suppose that ${\widehat K}_{l_1,l_2,\cdots,l_n}^U$ is not compact. 
Then, there exists 
$$b\in (\mathbb C^{l_1}\times\mathbb C^{l_2}\times\cdots\times\mathbb C^{l_n})
\setminus U$$ such that $$\inf_{a\in{\widehat K}_{l_1,l_2,\cdots,l_n}^U} |a-b| \, = \, 0$$ holds. 

\noindent 
{\sf Proof} : \ From Lemma 11, ${\widehat K}_{l_1,l_2,\cdots,l_n}^U$ is not a closed set 
of $\mathbb C^{l_1}\times\mathbb C^{l_2}\times\cdots\times\mathbb C^{l_n}$. 
So, there exist a sequence $\{a_m\}_{m=1}^\infty$ in ${\widehat K}_{l_1,l_2,\cdots,l_n}^U$ 
and $b\in (\mathbb C^{l_1}\times\mathbb C^{l_2}\times\cdots\times\mathbb C^{l_n})
\setminus {\widehat K}_{l_1,l_2,\cdots,l_n}^U$ 
such that $\lim_{ \, m\rightarrow\infty} \, a_m \, = \, b$ holds. 
Because ${\widehat K}_{l_1,l_2,\cdots,l_n}^U$ is a closed set of $U$,  
$b\not\in U$ holds. 
\hfill 
$\blacksquare$

\noindent 
{\bf Lemma 13} : 

Let $K$ be a compact subset of $U$. 
Let $$r \, :=\, 
\inf_{\, z\in K, \, w\in 
(\mathbb C^{l_1}\times\mathbb C^{l_2}\times\cdots\times\mathbb C^{l_n})\setminus U} 
\, |z-w|.$$ 
Then, for any  $a\in {\widehat K}_{l_1,l_2,\cdots,l_n}^U$ 
and $f\in O_{l_1,l_2,\cdots,l_n}(U)$, 
there exists $g\in O_{l_1,l_2,\cdots,l_n}
(\{ \, z\in \mathbb C^{l_1}\times\mathbb C^{l_2}\times\cdots\times\mathbb C^{l_n} 
\, | \, |z-a|<r \, \})$ such that for any multi-index $\alpha$, 
$\frac{\partial^{|\alpha|} f}{{\partial z}^\alpha}(a)
=\frac{\partial^{|\alpha|} g}{{\partial z}^\alpha}(a)$ holds. 

\noindent 
{\sf Proof} : \ Let $s\in (0,r)$. Then, from Cauchy inequality, 
there exists $c\in (0,+\infty)$ such that for any multi-index $\alpha$, 
$$\left( \, \left|\frac{\partial^{|\alpha|} f}{{\partial z}^\alpha}(a)\right|
\, \leq \, \right) \, \sup_{z\in K}\left|\frac{\partial^{|\alpha|} f}{{\partial z}^\alpha}(z)\right| 
\, \leq \, c\frac{\alpha !}{s^{|\alpha|}}$$
holds. Hence, $g: \, z \, \mapsto \, \sum_{\alpha}\frac{1}{\alpha !}
\frac{\partial^{|\alpha|} f}{{\partial z}^\alpha}(a)(z-a)^\alpha 
\, \in \, O_{l_1,l_2,\cdots,l_n}
(\{ \, z\in \mathbb C^{l_1}\times\mathbb C^{l_2}\times\cdots\times\mathbb C^{l_n} 
\, | \, |z-a|<r \, \})$ holds. 
\hfill 
$\blacksquare$

\noindent 
{\bf Lemma 14} : 

$\mathbb C^n$-holomorphy implies $\mathbb C^n$-holomorphic convexity. 

\noindent 
{\sf Proof} : \ Suppose that $U$ is not $\mathbb C^n$-holomorphically convex. 
Then, we show that $U$ is not a domain of $\mathbb C^n$-holomorphy. 
There exists a compact subset $K$ of $U$ such that 
${\widehat K}_{l_1,l_2,\cdots,l_n}^U$ is not compact. 
Let $$r \, :=\, 
\inf_{\, z\in K, \, w\in 
(\mathbb C^{l_1}\times\mathbb C^{l_2}\times\cdots\times\mathbb C^{l_n})\setminus U} 
\, |z-w|.$$ 
Then, from Lemma 12, there exist 
$a\in{\widehat K}_{l_1,l_2,\cdots,l_n}^U$ and  
$b\in (\mathbb C^{l_1}\times\mathbb C^{l_2}\times\cdots\times\mathbb C^{l_n})
\setminus U$ such that $$ |a-b| \, < \, \frac{r}{2}$$ holds. 
Hence, from Lemma 13 and Proposition 4, 
$U$ is not a domain of $\mathbb C^n$-holomorphy. 
\hfill 
$\blacksquare$

\noindent 
{\bf Lemma 15} : 

Let $\{K_m\}_{m=0}^\infty$ be a sequence 
of compact subsets of $U$. 
Let $\{p_m\}_{m=1}^\infty$ be a sequence in $U$. 
Suppose that $U\, = \, \cup_{m=0}^\infty \, ({K_m}^\circ)$ 
holds and for any nonnegative integer $m$, 
$K_m \, \subset \, K_{m+1}$ and $p_{m+1} 
\, \in \, K_{m+1}\setminus {\widehat{K_m}}_{l_1,l_2,\cdots,l_n}^U$ hold. 
Then, there exists $f\in O_{l_1,l_2,\cdots,l_n}(U)$ such that 
for any $m\in \mathbb N$, $m\leq |f(p_m)|$ holds. 

\noindent 
{\sf Proof} : \ From $p_{1} 
\, \not\in \, {\widehat{K_0}}_{l_1,l_2,\cdots,l_n}^U$, 
there exists $g_1 \, \in \, O_{l_1,l_2,\cdots,l_n}(U)$ such that 
$\sup_{\, w\in K_0} \, |g_1(w)| \, < \, |g_1(p_1)|$ holds. 
There exists $c_1\in (0,+\infty)$ such that 
$\sup_{\, w\in K_0} \, |c_1g_1(w)| \, < \, 1 \, < \, |c_1g_1(p_1)|$ holds. 
Then, there exists $k_1\in \mathbb N$ such that 
$\sup_{\, w\in K_0} \, |(c_1g_1(w))^{k_1}| \, \leq \, \frac{1}{2^{0}}$ 
and $2+\sum_{j=1}^{0}|(c_jg_j(p_1))^{k_j}| \, (=2) \, \leq \, |(c_1g_1(p_1))^{k_1}|$ hold. 
From $p_{2} 
\, \not\in \, {\widehat{K_1}}_{l_1,l_2,\cdots,l_n}^U$, 
there exists $g_2 \, \in \, O_{l_1,l_2,\cdots,l_n}(U)$ such that 
$\sup_{\, w\in K_1} \, |g_2(w)| \, < \, |g_2(p_2)|$ holds. 
There exists $c_2\in (0,+\infty)$ such that 
$\sup_{\, w\in K_1} \, |c_2g_2(w)| \, < \, 1 \, < \, |c_2g_2(p_2)|$ holds. 
Then, there exists $k_2\in \mathbb N$ such that 
$\sup_{\, w\in K_1} \, |(c_2g_2(w))^{k_2}| \, \leq \, \frac{1}{2^{1}}$ 
and $3+\sum_{j=1}^{1}|(c_jg_j(p_2))^{k_j}| \leq \, |(c_2g_2(p_2))^{k_2}|$ hold. 
Hereinafter, in the same manner, there exists a sequence
$\{(g_m,c_m,k_m)\}_{m=1}^\infty$
such that for any $m\in \mathbb N$, 
$g_m\in O_{l_1,l_2,\cdots,l_n}(U)$, 
$c_m\in (0,+\infty)$, 
$k_m\in \mathbb N$, 
$\sup_{w\in K_{m-1}}  |(c_mg_m(w))^{k_m}| \leq \frac{1}{2^{m-1}}$  
and $1+m+\sum_{j=1}^{m-1}|(c_jg_j(p_m))^{k_j}| \leq \, |(c_mg_m(p_m))^{k_m}|$ hold. 

For any $m\in \mathbb N$, 
$\sup_{w\in K_{m-1}}(\sum_{j=m}^\infty|(c_jg_j(w))^{k_j}|)
\, \leq \, \sum_{j=m}^\infty(\sup_{w\in K_{j-1}}|(c_jg_j(w))^{k_j}|) 
\, \leq \, \sum_{j=m}^\infty\frac{1}{2^{j-1}} \, = \, \frac{1}{2^{m-2}}$ 
holds. So, $f \, := \, \sum_{m=1}^\infty((c_mg_m)^{k_m}) 
\, \in \, O_{l_1,l_2,\cdots,l_n}(U)$ holds. 
For any $m\in \mathbb N$, 
$$1+m+|(c_mg_m(p_m))^{k_m}|$$
$$= \, 1+m+\left|f(p_m)-\left(\left(\sum_{j=1}^{m-1}((c_jg_j(p_m))^{k_j})\right)
+\left(\sum_{j=m+1}^{\infty}((c_jg_j(p_m))^{k_j})\right)\right)\right|$$
$$\leq \, 
1+m+|f(p_m)|+\left(\sum_{j=1}^{m-1} |(c_jg_j(p_m))^{k_j}|\right)
+\left(\sum_{j=m+1}^{\infty} |(c_jg_j(p_m))^{k_j}|\right)$$
$$\leq \, \left(\sum_{j=m+1}^{\infty} |(c_jg_j(p_m))^{k_j}|\right)
+|f(p_m)|+|(c_mg_m(p_m))^{k_m}|$$
and, so,   
$$1+m$$
$$\leq \, \left(\sum_{j=m+1}^{\infty} |(c_jg_j(p_m))^{k_j}|\right)+|f(p_m)|$$
$$\leq \, \left(\sum_{j=m+1}^{\infty} (\sup_{w\in K_{j-1}}|(c_jg_j(w))^{k_j}|)\right)+|f(p_m)|$$
$$\leq \, \left(\sum_{j=m+1}^{\infty} \frac{1}{2^{j-1}}\right)+|f(p_m)|$$
$$= \, \frac{1}{2^{m-1}}+|f(p_m)|$$
$$\leq \, 1+|f(p_m)|$$
hold. 
\hfill 
$\blacksquare$

\noindent 
{\bf Lemma 16} : 

Suppose that $U$ is $\mathbb C^n$-holomorphically convex. Suppose 
$U\, \not=\, \mathbb C^{l_1}\times\mathbb C^{l_2}\times\cdots\times\mathbb C^{l_n}$. 
Let $\{a_k\}_{k=1}^\infty$ be a sequence in $U$. 
For $k\in \mathbb N$, let
$$B_k \, \, \, := \, \, \, \{ \, z\in U \, | \, |a_k-z| \, < \, 
\inf_{w\in (\mathbb C^{l_1}\times\mathbb C^{l_2}\times\cdots\times\mathbb C^{l_n}) 
\setminus U} |a_k-w| \, \}.$$
Then, there exists $f\in O_{l_1,l_2,\cdots,l_n}(U)$ such that 
for any $k\in \mathbb N$, 
$$\sup_{z\in B_k} |f(z)| \, = \, +\infty$$
holds. 

\noindent 
{\sf Proof} : \ 
Let $$((q_1),(q_2,q_3),(q_4,q_5,q_6),(q_7,q_8,q_9,q_{10}),\cdots)$$
$$:= \, \, \, ((a_1),(a_1,a_2),(a_1,a_2,a_3),(a_1,a_2,a_3,a_4),\cdots).$$ 
Then, $\{q_m\}_{m=1}^\infty$ is a sequence in $U$ and 
for any $k\in \mathbb N$ and $l\in \mathbb N$, 
there exists $m \in \mathbb N$ 
such that $a_k=q_m$ and $l\leq m$ hold.  
Let $r_0:=1, R_0:=1$ and 
$$K_0$$
$$:= \, \, \, ( \, \cap_{\, w \, \in \, 
(\mathbb C^{l_1}\times\mathbb C^{l_2}\times\cdots\times\mathbb C^{l_n}) 
\setminus U} \, \{ \, z\in 
\mathbb C^{l_1}\times\mathbb C^{l_2}\times\cdots\times\mathbb C^{l_n} 
\, | \, r_0\leq |z-w| \, \} \, )$$
$$\cap \, \, \, \{ \, z\in 
\mathbb C^{l_1}\times\mathbb C^{l_2}\times\cdots\times\mathbb C^{l_n}
 \, | \, |z|\leq R_0 \, \}.$$  
Then, $K_0$ is a compact subset of $U$
and, so, ${\widehat{K_0}}_{l_1,l_2,\cdots,l_n}^U$ is a compact subset of $U$. 
Hence, there exists 
$p_1\, \in \, U\setminus {\widehat{K_0}}_{l_1,l_2,\cdots,l_n}^U$ 
such that 
$|q_1-p_1| \, < \, \inf_{\, w \, \in \, 
(\mathbb C^{l_1}\times\mathbb C^{l_2}\times\cdots\times\mathbb C^{l_n}) 
\setminus U} \, |q_1-w|$ and $\inf_{\, w \, \in \, 
(\mathbb C^{l_1}\times\mathbb C^{l_2}\times\cdots\times\mathbb C^{l_n}) 
\setminus U} \, |p_1-w| \, \leq \, \frac{1}{2}r_0$ hold. 
Let $r_1 \, := \, \inf_{\, w \, \in \, 
(\mathbb C^{l_1}\times\mathbb C^{l_2}\times\cdots\times\mathbb C^{l_n}) 
\setminus U} \, |p_1-w|, \, R_1\, :=\, \max\{|p_1|,2R_0\}$ and 
$$K_1$$
$$:= \, \, \, ( \, \cap_{\, w \, \in \, 
(\mathbb C^{l_1}\times\mathbb C^{l_2}\times\cdots\times\mathbb C^{l_n}) 
\setminus U} \, \{ \, z\in 
\mathbb C^{l_1}\times\mathbb C^{l_2}\times\cdots\times\mathbb C^{l_n} 
\, | \, r_1\leq |z-w| \, \} \, )$$
$$\cap \, \, \, \{ \, z\in 
\mathbb C^{l_1}\times\mathbb C^{l_2}\times\cdots\times\mathbb C^{l_n}
 \, | \, |z|\leq R_1 \, \}.$$ 
Then, $p_1\, \in \, K_1\setminus {\widehat{K_0}}_{l_1,l_2,\cdots,l_n}^U$, 
$0<r_1\leq\frac{1}{2}r_0<+\infty$ and $0<2R_0\leq R_1<+\infty$ 
hold. So, $K_1$ and ${\widehat{K_1}}_{l_1,l_2,\cdots,l_n}^U$ 
are compact subsets of $U$. 
Hence, there exists 
$p_2\, \in \, U\setminus {\widehat{K_1}}_{l_1,l_2,\cdots,l_n}^U$ 
such that 
$|q_2-p_2| \, < \, \inf_{\, w \, \in \, 
(\mathbb C^{l_1}\times\mathbb C^{l_2}\times\cdots\times\mathbb C^{l_n}) 
\setminus U} \, |q_2-w|$ and $\inf_{\, w \, \in \, 
(\mathbb C^{l_1}\times\mathbb C^{l_2}\times\cdots\times\mathbb C^{l_n}) 
\setminus U} \, |p_2-w| \, \leq \, \frac{1}{2}r_1$ hold. 
Let $r_2 \, := \, \inf_{\, w \, \in \, 
(\mathbb C^{l_1}\times\mathbb C^{l_2}\times\cdots\times\mathbb C^{l_n}) 
\setminus U} \, |p_2-w|, \, R_2\, :=\, \max\{|p_2|,2R_1\}$ and 
$$K_2$$
$$:= \, \, \, ( \, \cap_{\, w \, \in \, 
(\mathbb C^{l_1}\times\mathbb C^{l_2}\times\cdots\times\mathbb C^{l_n}) 
\setminus U} \, \{ \, z\in 
\mathbb C^{l_1}\times\mathbb C^{l_2}\times\cdots\times\mathbb C^{l_n} 
\, | \, r_2\leq |z-w| \, \} \, )$$
$$\cap \, \, \, \{ \, z\in 
\mathbb C^{l_1}\times\mathbb C^{l_2}\times\cdots\times\mathbb C^{l_n}
 \, | \, |z|\leq R_2 \, \}.$$
Then, $p_2\, \in \, K_2\setminus {\widehat{K_1}}_{l_1,l_2,\cdots,l_n}^U$, 
$0<r_2\leq\frac{1}{2}r_1<+\infty$ and $0<2R_1\leq R_2<+\infty$ 
hold. Hereinafter, in the same manner, there exist sequences 
$\{(r_m,R_m,K_m)\}_{m=0}^\infty$ and $\{p_m\}_{m=1}^\infty$ such that 
for any nonnegative integer $m$, 
$0<r_{m+1}\leq\frac{1}{2}r_m<+\infty$, 
$0<2R_m\leq R_{m+1}<+\infty$,  
$$K_m$$
$$= \, \, \, ( \, \cap_{\, w \, \in \, 
(\mathbb C^{l_1}\times\mathbb C^{l_2}\times\cdots\times\mathbb C^{l_n}) 
\setminus U} \, \{ \, z\in 
\mathbb C^{l_1}\times\mathbb C^{l_2}\times\cdots\times\mathbb C^{l_n} 
\, | \, r_m\leq |z-w| \, \} \, )$$
$$\cap \, \, \, \{ \, z\in 
\mathbb C^{l_1}\times\mathbb C^{l_2}\times\cdots\times\mathbb C^{l_n}
 \, | \, |z|\leq R_m \, \},$$  
$p_{m+1}\, \in \, K_{m+1}\setminus {\widehat{K_m}}_{l_1,l_2,\cdots,l_n}^U$
and $$|q_{m+1}-p_{m+1}| \, < \, \inf_{\, w \, \in \, 
(\mathbb C^{l_1}\times\mathbb C^{l_2}\times\cdots\times\mathbb C^{l_n}) 
\setminus U} \, |q_{m+1}-w|$$ hold.  Then, from Lemma 15, 
there exists $f\in O_{l_1,l_2,\cdots,l_n}(U)$ such that 
for any $m\in\mathbb N$, $m\leq |f(p_m)|$ holds. 

Let $k\in \mathbb N$. We show $\sup_{z\in B_k} |f(z)| \, = \, +\infty$. 
Let $l\in \mathbb N$. Then, there exists $m \in \mathbb N$ 
such that $a_k=q_m$ and $l\leq m$ hold.  Hence, 
$|a_k-p_{m}| \, < \, \inf_{\, w \, \in \, 
(\mathbb C^{l_1}\times\mathbb C^{l_2}\times\cdots\times\mathbb C^{l_n}) 
\setminus U} \, |a_k-w|$ and $l\leq |f(p_m)|$ hold.  So, 
$l\leq \sup_{z\in B_k} |f(z)|$ holds. Therefore, 
$\sup_{z\in B_k} |f(z)| \, = \, +\infty$ holds. 
\hfill 
$\blacksquare$

\vspace*{0.8em}

\noindent 
{\sf Proof of Theorem 8} : \ 
Suppose that $U$ is $\mathbb C^n$-holomorphically convex. 
We show that  $U$ is a domain of $\mathbb C^n$-existence. 
When $U\, =\, \emptyset$ 
or $U\, =\, \mathbb C^{l_1}\times\mathbb C^{l_2}\times\cdots\times\mathbb C^{l_n}$ 
holds, it is obvious. 
Suppose $U\, \not=\, \emptyset$ 
and $U\, \not=\, \mathbb C^{l_1}\times\mathbb C^{l_2}\times\cdots\times\mathbb C^{l_n}$. 
Then, there exists a sequence $\{a_k\}_{k=1}^\infty$ in $U$ such that 
$$U \, = \, \overline{\{a_k\}_{k=1}^\infty}$$ 
holds. For $k\in \mathbb N$, let
$$B_k \, \, \, := \, \, \, \{ \, z\in U \, | \, |a_k-z| \, < \, 
\inf_{w\in (\mathbb C^{l_1}\times\mathbb C^{l_2}\times\cdots\times\mathbb C^{l_n}) 
\setminus U} |a_k-w| \, \}.$$
Then, from Lemma 16, there exists $f\in O_{l_1,l_2,\cdots,l_n}(U)$ such that 
for any $k\in \mathbb N$, 
$$\sup_{z\in B_k} |f(z)| \, = \, +\infty$$
holds. 

We show that $U$ is the domain of $\mathbb C^n$-existence of $f$. 
Suppose that $U$ is not the domain of $\mathbb C^n$-existence of $f$. 
Then, there exist open sets $V$ and $W$ of 
$\mathbb C^{l_1}\times\mathbb C^{l_2}\times\cdots\times\mathbb C^{l_n}$ 
and $g\in  O_{l_1,l_2,\cdots,l_n}(V)$ 
such that $V$ is connected and $\emptyset\not=V\setminus U$,  
$\emptyset\not=W\subset U\cap V$ and 
$f_{\upharpoonright W}=g_{\upharpoonright W}$ hold. 
Let $$\tilde{W} \, \, \, := \, \, \, \{\, w\in U\cap V\, |\, 
\exists \, r\in(0,+\infty), \, \forall \,  z\in U\cap V: \, 
[ \, |z-w|<r \, \Rightarrow \, f(z)=g(z) \, ]   \, \}.$$ 
So, $\emptyset\not=\tilde{W}\subsetneq V$ holds and 
$\tilde{W}$ is an open set of $V$. 
Because $V$ is connected, $\tilde{W}$ is not an closed set of $V$. 
Hence, there exists $b \, \in \, 
\left(V\cap \overline{\tilde{W}}\right)\setminus \tilde{W}$. 
We show $b\not\in U$. Suppose $b\in U$. 
Then, $b\in (U\cap V)\cap \overline{\tilde{W}}$ holds. 
Hence, from Proposition 4, $b\in \tilde{W}$ holds. 
It is a contradiction. So, $b\not\in U$ holds. 
Therefore, $$b \, \in \, 
\left(V\cap \overline{\tilde{W}}\right)\setminus U$$ 
holds. Let $\varepsilon\in (0,+\infty)$. 
Let $\delta \, := \, \min\{\varepsilon, 
\inf_{w\in (\mathbb C^{l_1}\times\mathbb C^{l_2}\times\cdots\times\mathbb C^{l_n}) 
\setminus V} |b-w| \}$. Then, there exists $a\in \tilde{W}$ 
such that $|a-b|<\frac{\delta}{4}$ holds. 
Further, there exists $k\in \mathbb N$ 
such that $|a_k-a|<\frac{\delta}{4}$ and 
$a_k\in \tilde{W}$ hold. For any $z\in B_k$,  
$|a_k-z| \, < \, 
\inf_{w\in (\mathbb C^{l_1}\times\mathbb C^{l_2}\times\cdots\times\mathbb C^{l_n}) 
\setminus U} |a_k-w| 
\, \leq \, |a_k-b| \, < \, \frac{\delta}{2}$ and, so,  
$|b-z| \, < \, \delta \, \leq \, 
\inf_{w\in (\mathbb C^{l_1}\times\mathbb C^{l_2}\times\cdots\times\mathbb C^{l_n}) 
\setminus V} |b-w|$ hold. Hence, $B_k\subset V$ holds. 
$B_k\subset U\cap V$ and $a_k\in B_k\cap \tilde{W}$ hold and 
$B_k$ is connected. So, from Proposition 4, $B_k\subset \tilde{W}$ holds. 
Hence, because $z\in B_k$ implies  
$|b-z| \, < \, \delta \, \leq \, \varepsilon$, 
$$( \, +\infty \, = \, ) \, \sup_{z\in B_k} |f(z)| \, = \, 
\sup_{z\in B_k} |g(z)| \, \leq \, 
\sup_{z \, \in \, \{ \, w\in V \, | \, |b-w|<\varepsilon \, \}}|g(z)|$$
holds.  Therefore, for any $\varepsilon\in (0,+\infty)$, 
$\sup_{\, z \, \in \, \{ \, w\in V \, | \, |b-w|<\varepsilon \, \}} \, |g(z)| 
\, = \, +\infty$ holds. 
However, since $b\in V$ and $g\in  O_{l_1,l_2,\cdots,l_n}(V)$ hold, it is a contradiction. 
So, $U$ is the domain of $\mathbb C^n$-existence of $f$. 

Because $\mathbb C^n$-holomorphic convexity 
implies $\mathbb C^n$-existence, from Lemmas 6 and 14,
it follows. 
\hfill 
$\blacksquare$


\vfill

Acknowledgment: 

As in Comment, Professor Naohiko Kasuya suggested it. 

This work was supported by JSPS KAKENHI Grant Number JP16K05245. 

\newpage 

{\bf References}

\vspace*{0.8em}

[1] H. Cartan and P. Thullen, Zur Theorie der Singularitaten der Funktionen mehrerer komplexen Veranderlichen (German), {\it Math. Ann.}, 106 (1932),
617-647.

[2] B. W. Glickfeld, The Riemann sphere of a commutative Banach algebra, 
{\it Trans. Amer. Math. Soc.}, 134 (1968), 1-28.

[3] L. Hormander, 
$L^2$ estimates and existence theorems for the $\overline \partial$ operator,
{\it Acta Math.}, 113 (1965), 89-152.

[4] S. Kobayashi, Manifolds over function algebras and mapping spaces,
{\it Tohoku Math. J.}, 41 (1989), 263-282.

[5] L. Lempert, The Dolbeault complex in infinite dimensions, 
{\it J. Amer. Math. Soc.}, 11 (1998), 485-520.

[6] E. R. Lorch, The theory of analytic functions in normed Abelian vector
rings, {\it Trans. Amer. Math. Soc.}, 54 (1943), 414-425.

[7] A. Mallios and E. E. Rosinger, Space-time foam dense singularities
and de Rham cohomology, {\it Acta Appl. Math.}, 67 (2001), 59-89.

[8] P. Manoharan, A nonlinear version of Swanfs theorem, 
{\it Math. Z.}, 209
(1992), 467-479.

[9] P. Manoharan, Generalized Swanfs theorem and its application, 
{\it Proc. Amer. Math. Soc.}, 123 (1995), 3219-3223.

[10] P. Manoharan, A characterization for spaces of sections, 
{\it Proc. Amer. Math. Soc.}, 126 (1998), 1205-1210.

[11] A. S. Morye, Note on the Serre-Swan theorem, 
{\it Math. Nachr.}, 286
(2013), 272-278.

[12] T. Ohsawa and K. Takegoshi, On the extension 
of $L^2$ holomorphic
functions, {\it Math. Z.}, 195 (1987), 197-204.

[13] M. H. Papatriantafillou, Partitions of unity on A-manifolds, {\it Internat.
J. Math.}, 9 (1998), 877-883.

[14] R. G. Swan, Vector bundles and projective modules, {\it Trans. Amer.
Math. Soc.}, 105 (1962), 264-277.

[15] H. Yagisita, Finite-dimensional complex manifolds on commutative Banach algebras and continuous families of compact complex manifolds, 
{\it Complex Manifolds}, 6 (2019), 228-264. 

[16] H. Yagisita, Holomorphic differential forms of complex manifolds on
commutative Banach algebras and a few related problems, {\it arXiv.org}.







\end{document}